\documentclass[12pt]{article}
\usepackage{amsmath}
\usepackage{amsfonts}
\usepackage{amssymb}
\usepackage[doublespacing]{setspace}

\begin{document}

\begin{center}
{\LARGE Oscillation of a Class of Impulsive Differential Equations with
Continuous and Piecewise Constant Arguments}

\bigskip

Fatma KARAKO\c{C}$^{1}\footnotetext[1]{%
Corresponding Author: fkarakoc@ankara.edu.tr
\par
Phone: +90(312)2126720
\par
Fax: +90(312)2235000
\par
{}
\par
{}}$

Department of Mathematics, Faculty of Sciences, Ankara University, Turkey
\end{center}

\noindent \textbf{Abstract }{\small A class of first order linear impulsive
differential equation with continuous and piecewise constant arguments is
studied. Sufficient conditions for the oscillation of the solutions are
obtained.\medskip }



\section*{\textbf{1. Introduction}}

\noindent In this paper, we consider an impulsive differential equation with
continuous and piecewise constant arguments of the form%
\begin{equation}
x^{\prime }\left( t\right) +a\left( t\right) x\left( t\right) +b(t)x(t-\tau
)+c(t)x([t-1])=0,\ t\neq t_{i},~t\geq t_{0}>0,  \label{1}
\end{equation}%
\begin{equation}
\Delta x\left( t_{i}\right) =b_{i}x\left( t_{i}\right) ,\ i=1,2,...,\
\label{2}
\end{equation}%
where $a\in C([0,\infty ),%
\mathbb{R}
),~b,c\in C([0,\infty ),[0,\infty ))$, $\tau \in
\mathbb{R}
^{+}$ is a fixed constant, $\left[ .\right] $ denotes the greatest integer
function, $\{t_{i}\}$ is a sequence of real numbers such that $%
~0<t_{0}<t_{1}<t_{2}<...<t_{j}<t_{j+1}<...,$ and $\underset{i\rightarrow
\infty }{\lim }t_{i}=\infty ,~\Delta x\left( t_{i}\right) =x\left(
t_{i}^{+}\right) -x\left( t_{i}^{-}\right) ,\ x\left( t_{i}^{+}\right) =\lim
\limits_{t\rightarrow t_{i}^{+}}x\left( t\right) ,\ x\left( t_{i}^{-}\right)
=\lim \limits_{t\rightarrow t_{i}^{-}}x\left( t\right) ,$ $b_{i}\neq 1$, $%
i=1,2,...,$ are constants.

\noindent Since 1980's differential equations with piecewise constant
arguments have been attracted great deal of attention of researchers in
mathematical and some of the others fields in science. Piecewise constant
systems exist in a widely expanded areas such as biomedicine, chemistry,
mechanical engineering, physics, etc. These kind of equations such as Eq.(%
\ref{1}) are similar in structure to those found in certain
sequential-continuous models of disease dynamics \cite{1}. In 1994, Dai and
Sing \cite{2} studied the oscillatory motion of spring-mass systems with
subject to piecewise constant forces of the form $f(x[t])$ or $f([t]).$
Later, they improved an analytical and numerical method for solving linear
and nonlinear vibration problems and they showed that a function $%
f([N(t)]/N) $ is a good approximation to the given continuous function $f(t)$
if $N$ is sufficiently large \cite{3}.

\noindent In 1984, Cooke and Wiener \cite{4} studied oscillatory and
periodic solutions of a linear differential equation with piecewise constant
argument and they note that such equations are comprehensively related to
impulsive and difference equations. After this work, oscillatory and
periodic solutions of linear differential equations with piecewise constant
arguments have been dealt with by many authors \cite{5,6,7} and the
references cited therein.

\noindent On the other hand, in 1994, the case of studying discontinuous
solutions of differential equations with piecewise continuous arguments has
been proposed as an open problem by Wiener \cite{8}. Due to this open
problem, some impulsive differential equations with piecewise constant
arguments have been studied \cite{9,10,11}. Moreover, the monographs \cite%
{12,13} includes many results on the theory of differential equations with
piecewise constant arguments.

\noindent Now, our aim is to consider the Wiener's open problem for the
equation (\ref{1})-(\ref{2}). Moreover, as we know there is only one work on
nonimpulsive delay differential equations with continuous and piecewise
constant arguments \cite{14}. In this respect, we obtain sufficient
conditions for the existence of oscillatory solutions of Eq. (\ref{1})-(\ref%
{2}).\medskip \

\noindent \textbf{Definition 1 }It is said that a function $x$ defined on
the set $\{-1\} \cup \lbrack -\tau ,\infty )$ is a solution of Eq. (\ref{1}%
)-(\ref{2}) if it satisfies the following conditions:

\noindent $(D_{1})$ $x(t)$ is continuous on $[-\tau ,\infty )~$with the
possible exception of the points $t_{i},~i=1,2,...$\newline
$(D_{2})~x(t)$ is right continuous and has left-hand limit at the points $%
t_{i},~i=1,2,...$\newline
$(D_{3})$ $x(t)$ differentiable and satisfies (\ref{1}) for any $t\in
\mathbb{R}
^{+},$ with the possible\ exception of the points $t_{i},~i=1,2,...,$\ and $%
[t]\in \lbrack 0,\infty ),$ where one-sided derivatives exist,\newline
$(D_{4})$ $x(t)$ satisfies (\ref{2}) at the points $t_{i},~i=1,2,...\medskip
$

\noindent \textbf{Definition 2. }A function $x\left( t\right) $ is called
oscillatory if it is neither positive nor negative for $t\geq T$ where $T$
is sufficiently large. Otherwise, the solution is called
nonoscillatory.\medskip

\noindent \textbf{Remark 1.} In this paper we assume that $-\infty <b_{i}<1$
for all $i=1,2,...$ Otherwise, from the impulse conditions (\ref{2}) it is
obtained that the solutions are already oscillatory.

\noindent \textbf{Remark 2.} We assume that $b(t)\not \equiv 0$ or $c(t)\not
\equiv 0.$ If $b(t)\equiv 0$ and $c(t)\equiv 0,$ then Eq. (\ref{1})-(\ref{2}%
) reduces an ordinary differential equation with impulses. The results on
the oscillation of impulsive ordinary differential equations can be found in
the survey paper \cite{15}.

\section*{\textbf{2. Main Results}}

\noindent In this paper we also consider following differential inequalities.%
\begin{eqnarray}
x^{\prime }\left( t\right) +a\left( t\right) x\left( t\right) +b(t)x(t-\tau
)+c(t)x([t-1]) &\leq &0,\ t\neq t_{i},~t\geq t_{0}>0,  \label{3} \\
\Delta x\left( t_{i}\right) &=&b_{i}x\left( t_{i}\right) ,\ i=1,2,...,
\notag
\end{eqnarray}%
and%
\begin{eqnarray}
x^{\prime }\left( t\right) +a\left( t\right) x\left( t\right) +b(t)x(t-\tau
)+c(t)x([t-1]) &\geq &0,\ t\neq t_{i},~t\geq t_{0}>0,  \label{4} \\
\Delta x\left( t_{i}\right) &=&b_{i}x\left( t_{i}\right) ,\ i=1,2,...  \notag
\end{eqnarray}%
The main tools for the proofs of our results are following differential
equation and inequalities.%
\begin{equation}
y^{\prime }\left( t\right) +a\left( t\right) y\left( t\right) +B(t)y(t-\tau
)+C(t)y([t-1])=0,~t\geq t_{0}+\max \{ \tau ,2\}  \label{5}
\end{equation}%
\begin{equation}
y^{\prime }\left( t\right) +a\left( t\right) y\left( t\right) +B(t)y(t-\tau
)+C(t)y([t-1])\leq 0,~t\geq t_{0}+\max \{ \tau ,2\},  \label{6yy}
\end{equation}%
\begin{equation}
y^{\prime }\left( t\right) +a\left( t\right) y\left( t\right) +B(t)y(t-\tau
)+C(t)y([t-1])\geq 0,~t\geq t_{0}+\max \{ \tau ,2\},  \label{7y}
\end{equation}%
where
\begin{equation}
B(t)=\prod \limits_{t-\tau <t_{j}\leq t}(1-b_{j})b(t),~t\geq t_{0}+\max \{
\tau ,2\},  \label{8y}
\end{equation}%
and%
\begin{equation}
C(t)=\prod \limits_{[t-1]<t_{j}\leq t}(1-b_{j})c(t),~t\geq t_{0}+\max \{
\tau ,2\}.  \label{9y}
\end{equation}%
The following theorem is a generalization of Theorem 1 in \cite{16} to
impulsive delay differential equations with continuous and piecewise
constant arguments.\medskip

\noindent \textbf{Theorem 1.} $(i)$~Inequality (\ref{3}) has no eventually
positive solution if and only if inequality (\ref{6yy}) has no eventually
positive solution.

\noindent $(ii)$ Inequality (\ref{4}) has no eventually negative solution if
and only if inequality (\ref{7y}) has no eventually negative solution.

\noindent $(iii)$ All solutions of the equation (\ref{1})-(\ref{2}) are
oscillatory if and only if all solutions of equation (\ref{5}) are
oscillatory.\medskip

\noindent \textbf{Proof. }We will prove $(i)$ since\ the proofs of $(ii)$
and $(iii)$ are similar to proof of $(i).$ Let $x(t)$ be an eventually
positive solution of inequality (\ref{3}) such that $x(t)>0,~x(t-\tau
)>0,~x([t-1])>0$ for $t>T\geq t_{0}+\max \{ \tau ,2\},$ where $T$ is
sufficiently large. Set $y(t)=\prod \limits_{T<t_{j}\leq t}(1-b_{j})x(t).$
Since $1-b_{j}>0,$ it is clear that $y(t)>0,~y(t-\tau )>0,$ and $y([t-1])>0$
for $t>T.$ Now, we will show that $y(t)$ is a solution of inequality (\ref%
{6yy}). From (\ref{8y}), (\ref{9y}), and (\ref{3}) we obtain that

$y^{\prime }\left( t\right) +a\left( t\right) y\left( t\right) +B(t)y(t-\tau
)+C(t)y([t-1])~\  \  \  \  \  \  \  \  \  \  \  \  \  \  \  \  \  \  \  \  \  \  \  \  \  \  \  \  \  \  \
\  \  \  \  \  \  \  \  \  \  \  \  \ $%
\begin{eqnarray*}
&=&\prod \limits_{T<t_{j}\leq t}(1-b_{j})x^{\prime }(t)+a(t)\prod
\limits_{T<t_{j}\leq t}(1-b_{j})x(t) \\
&&+\prod \limits_{t-\tau <t_{j}\leq t}(1-b_{j})b(t)\prod
\limits_{T<t_{j}\leq t-\tau }(1-b_{j})x(t-\tau ) \\
&&+\prod \limits_{[t-1]<t_{j}\leq t}(1-b_{j})c(t)\prod \limits_{T<t_{j}\leq
\lbrack t-1]}(1-b_{j})x([t-1]) \\
&\leq &\prod \limits_{T<t_{j}\leq t}(1-b_{j})\left[ x^{\prime }\left(
t\right) +a\left( t\right) x\left( t\right) +b(t)x(t-\tau )+c(t)x([t-1])%
\right] \\
&\leq &0.
\end{eqnarray*}%
So, $y(t)$ is an eventually positive solution of inequality (\ref{6yy}). On
the other hand, from (\ref{2}), we have%
\begin{eqnarray*}
y(t_{i}^{-}) &=&\prod \limits_{T<t_{j}\leq t_{i-1}}(1-b_{j})x(t_{i}^{-}) \\
&=&\prod \limits_{T<t_{j}\leq t_{i}}(1-b_{j})x(t_{i}) \\
&=&y(t_{i})
\end{eqnarray*}%
and%
\begin{eqnarray*}
y(t_{i}^{+}) &=&\prod \limits_{T<t_{j}\leq t_{i}}(1-b_{j})x(t_{i}^{+}) \\
&=&\prod \limits_{T<t_{j}\leq t_{i}}(1-b_{j})x(t_{i}) \\
&=&y(t_{i}).
\end{eqnarray*}%
So, $y(t)$ is continuous at the impulse points.

\noindent Now, let $y(t)$ be an eventually positive solution of inequality (%
\ref{6yy}). Then $y(t)>0,~y(t-\tau )>0,$ and $y([t-1])>0$ for $t>T.$ We will
show that $x(t)$ is an eventually positive solution of inequality (\ref{3}).
From (\ref{8y}), (\ref{9y}), and (\ref{6yy}) we obtain that

$x^{\prime }\left( t\right) +a\left( t\right) x\left( t\right) +b(t)x(t-\tau
)+c(t)x([t-1])~\  \  \  \  \  \  \  \  \  \  \  \  \  \  \  \  \  \  \  \  \  \  \  \  \  \  \  \  \  \  \
\  \  \  \  \  \  \  \  \  \  \  \  \  \  \  \  \  \  \  \  \  \  \ $%
\begin{eqnarray*}
&=&\prod \limits_{T<t_{j}\leq t}(1-b_{j})^{-1}y^{\prime }(t)+a(t)\prod
\limits_{T<t_{j}\leq t}(1-b_{j})^{-1}y(t) \\
&&+b(t)\prod \limits_{T<t_{j}\leq t-\tau }(1-b_{j})^{-1}y(t-\tau )+c(t)\prod
\limits_{T<t_{j}\leq \lbrack t-1]}(1-b_{j})^{-1}y([t-1]) \\
&=&\prod \limits_{T<t_{j}\leq t}(1-b_{j})^{-1}\left[ y^{\prime }\left(
t\right) +a\left( t\right) y\left( t\right) +B(t)y(t-\tau )+C(t)y([t-1])%
\right] \\
&\leq &0.
\end{eqnarray*}%
Moreover,%
\begin{eqnarray*}
x(t_{i}^{-}) &=&\prod \limits_{T<t_{j}\leq t_{i-1}}(1-b_{j})^{-1}y(t_{i}^{-})
\\
&=&\prod \limits_{T<t_{j}\leq t_{i}}(1-b_{j})^{-1}(1-b_{i})y(t_{i}) \\
&=&(1-b_{i})x(t_{i})
\end{eqnarray*}%
and%
\begin{equation*}
x(t_{i}^{+})=\prod \limits_{T<t_{j}\leq
t_{i}}(1-b_{j})^{-1}y(t_{i}^{+})=x(t_{i}).
\end{equation*}%
So, $x(t)$ is an eventually positive solution of inequality (\ref{3}). The
proof is complete.\medskip

\noindent Following we give several sufficient conditions for the
oscillation of equation (\ref{1})-(\ref{2}).\medskip

\noindent \textbf{Theorem 2. }If one of the following conditions is
satisfied then every solution of equation (\ref{1})-(\ref{2}) is oscillatory:

\begin{equation}
\lim \limits_{t\rightarrow \infty }\sup ~\int \limits_{t-l}^{t}\left(
\prod \limits_{s-\tau <t_{j}\leq s}(1-b_{j})\right) b(s)\exp \left(
~\int \limits_{s-\tau }^{s}a(u)du\right) ds>1,  \label{10y}
\end{equation}%
\begin{equation}
\lim \limits_{n\rightarrow \infty }\sup ~\int \limits_{n+1-l}^{n+1}\left(
\prod \limits_{n-1<t_{j}\leq s}(1-b_{j})\right) c(s)\exp \left(
~\int \limits_{n-1}^{s}a(u)du\right) ds>1,  \label{11y}
\end{equation}%
where $l=\min \{ \tau ,1\}.$

\noindent \textbf{Proof. }Let conditions (\ref{10y}) or (\ref{11y})\textbf{\
}is satisfied. We shall prove that the existence of eventually positive (or
negative) solutions leads to a contradiction. Let $x(t)$ be an eventually
positive solution of equation (\ref{1})-(\ref{2}). Then $y(t)=\prod%
\limits_{T<t_{j}\leq t}(1-b_{j})x(t)$ is an eventually positive solution of
equation (\ref{5}) such that $y(t)>0,~y(t-\tau )>0,~y([t-1])>0$ for $%
n+1>t\geq n>T.$ Taking%
\begin{equation}
z(t)=y(t)\exp \left( ~\int \limits_{T}^{t}a(s)ds\right) ,~t>T,  \label{6y}
\end{equation}%
it is obtained from equation (\ref{5}) that%
\begin{equation}
z^{\prime }(t)=-\left \{ B(t)z(t-\tau )\exp \left( ~\int \limits_{t-\tau
}^{t}a(s)ds\right) +C(t)z([t-1])\exp \left( ~\int
\limits_{[t-1]}^{t}a(s)ds\right) \right \}  \label{7}
\end{equation}%
for $n+1>t\geq n>T.$ Since $B(t),~C(t)\geq 0$ for $t\in
\mathbb{R}
$ and $z(t-\tau ),~z([t-1])\geq 0$ for $n+1>t\geq n>T,~$\ we get $z(t)$ is
nonincreasing for $t>T.$

\noindent Now, we consider two cases:

\textit{Case 1.} $\tau >1.$ Then it is clear that $z(t-\tau )\geq z(t-1)$
and $z([t-1])\geq z(t-1)$ for $t>T.$

\noindent Using (\ref{7}), we obtain that%
\begin{eqnarray}
0 &=&z^{\prime }(t)+B(t)z(t-\tau )\exp \left( ~\int \limits_{t-\tau
}^{t}a(s)ds\right) +C(t)z([t-1])\exp \left( ~\int
\limits_{[t-1]}^{t}a(s)ds\right)  \notag \\
&\geq &z^{\prime }(t)+z(t-1)P(t),  \label{8}
\end{eqnarray}%
where%
\begin{equation}
P(t)=B(t)\exp \left( ~\int \limits_{t-\tau }^{t}a(s)ds\right) +C(t)\exp
\left( ~\int \limits_{[t-1]}^{t}a(s)ds\right) .  \label{6}
\end{equation}%
Integrating inequality (\ref{8}) from $t-1$ to $t$, we get%
\begin{equation*}
z(t)-z(t-1)+\int \limits_{t-1}^{t}P(s)z(s-1)ds\leq 0.
\end{equation*}%
Since $z(t)$ is nonincreasing for $t>T,$ from the above inequality, we
obtain that%
\begin{equation*}
z(t)+z(t-1)\left[ \int \limits_{t-1}^{t}P(s)ds-1\right] \leq 0
\end{equation*}%
and so, we have%
\begin{equation*}
\int \limits_{t-1}^{t}P(s)ds\leq 1.
\end{equation*}%
Using (\ref{6}), (\ref{8y}), and (\ref{9y}), we obtain from the above
inequality that%
\begin{equation}
\int \limits_{t-1}^{t}\prod \limits_{s-\tau <t_{j}\leq s}(1-b_{j})b(s)\exp
\left( ~\int \limits_{s-\tau }^{s}a(u)du\right) ds\leq 1,  \label{13y}
\end{equation}%
and%
\begin{equation*}
\int \limits_{t-1}^{t}\prod \limits_{[s-1]<t_{j}\leq s}(1-b_{j})c(s)\exp
\left( ~\int \limits_{[s-1]}^{s}a(u)du\right) ds\leq 1.
\end{equation*}%
It is clear that inequality (\ref{13y}) contradicts (\ref{10y}). On the
other hand, integrating inequality (\ref{8}) from $n$ to $n+1$, we get%
\begin{equation*}
z(n+1)-z(n)+\int \limits_{n}^{n+1}P(s)z(s-1)ds\leq 0.
\end{equation*}%
Since $z(t)$ is nonincreasing for $t>T,$ from the above inequality, we
obtain that%
\begin{equation*}
z(n+1)+z(n)\left[ \int \limits_{n}^{n+1}P(s)ds-1\right] \leq 0
\end{equation*}%
and so, we have
\begin{equation*}
\int \limits_{n}^{n+1}P(s)ds\leq 1.
\end{equation*}%
In view of (\ref{6}), (\ref{8y}), and (\ref{9y}), we obtain from the above
inequality that%
\begin{equation*}
\int \limits_{n}^{n+1}\prod \limits_{s-\tau <t_{j}\leq s}(1-b_{j})b(s)\exp
\left( ~\int \limits_{s-\tau }^{s}a(u)du\right) ds\leq 1,
\end{equation*}%
and%
\begin{equation}
\int \limits_{n}^{n+1}\prod \limits_{[s-1]<t_{j}\leq s}(1-b_{j})c(s)\exp
\left( ~\int \limits_{[s-1]}^{s}a(u)du\right) ds\leq 1.  \label{14y}
\end{equation}%
Since $n\leq s<n+1,~$(\ref{14y}) contradicts (\ref{11y}).

\textit{Case 2}. $\tau \leq 1.$ Then $z(t-\tau )\leq z([t-1])$ for $%
n+1>t\geq n>T,$ and from (\ref{7}), we obtain that%
\begin{eqnarray}
0 &=&z^{\prime }(t)+B(t)z(t-\tau )\exp \left( ~\int \limits_{t-\tau
}^{t}a(s)ds\right) +C(t)z([t-1])\exp \left( ~\int
\limits_{[t-1]}^{t}a(s)ds\right)  \notag \\
&\geq &z^{\prime }(t)+z(t-\tau )P(t),  \label{9}
\end{eqnarray}%
where $P(t)$ is defined in (\ref{6}). Integrating inequality (\ref{9}) from $%
t-\tau $ to $t$, we get%
\begin{equation*}
z(t)-z(t-\tau )+\int \limits_{t-\tau }^{t}P(s)z(s-\tau )ds\leq 0.
\end{equation*}%
Since $z(t)$ is nonincreasing for $t>T,$ from the above inequality, we
obtain that
\begin{equation*}
z(t)+z(t-\tau )\left[ \int \limits_{t-\tau }^{t}P(s)ds-1\right] \leq 0
\end{equation*}%
and so, we have%
\begin{equation*}
\int \limits_{t-\tau }^{t}P(s)ds\leq 1.
\end{equation*}%
Using (\ref{6}), (\ref{8y}), and (\ref{9y}), we obtain from the above
inequality that%
\begin{equation*}
\int \limits_{t-\tau }^{t}\left( \prod \limits_{s-\tau <t_{j}\leq
s}(1-b_{j})\right) b(s)\exp \left( ~\int \limits_{s-\tau }^{s}a(u)du\right)
ds\leq 1
\end{equation*}%
which contradicts (\ref{10y}). On the other hand, integrating inequality (%
\ref{9}) from $n+1-\tau $ to $n+1$, we get%
\begin{equation*}
z(n+1)-z(n+1-\tau )+\int \limits_{n+1-\tau }^{n+1}P(s)z(s-\tau )ds\leq 0.
\end{equation*}%
Since $z(t)$ is nonincreasing for $t>T,$ from the above inequality, we
obtain that%
\begin{equation*}
\int \limits_{n+1-\tau }^{n+1}P(s)ds\leq 1.
\end{equation*}%
In view of (\ref{6}), (\ref{8y}), and (\ref{9y}), we obtain from the above
inequality that%
\begin{equation*}
\int \limits_{n+1-\tau }^{n+1}\left( \prod \limits_{n-1<t_{j}\leq
s}(1-b_{j})\right) c(s)\exp \left( ~\int \limits_{n-1}^{s}a(u)du\right)
ds\leq 1
\end{equation*}%
which contradicts (\ref{11y}).

\noindent If $x(t)\ $is an eventually negative solution of equation (\ref{1}%
)-(\ref{2}), then $-x(t)$ is an eventually positive solution of equation (%
\ref{1})-(\ref{2}) and we obtain same contradiction. So, the proof is
complete.\medskip

\noindent \textbf{Corollary 1.} Assume that $b(t)\neq 0$, $c(t)\equiv 0$ and
that%
\begin{equation}
\lim \limits_{t\rightarrow \infty }\sup ~\int \limits_{t-\tau }^{t}\left(
\prod \limits_{s-\tau <t_{j}\leq s}(1-b_{j})\right) b(s)\exp \left( ~\int
\limits_{s-\tau }^{s}a(u)du\right) ds>1.  \label{10}
\end{equation}%
Then every solution of Eq. (\ref{1})-(\ref{2}) is oscillatory.

\noindent \textbf{Remark 3.} If $b(t)\neq 0$ and $c(t)\equiv 0,$ then Eq. (%
\ref{1})-(\ref{2}) reduces to a delay differential equation with impulses.
Condition (\ref{10}) is similar to hypothesis of Theorem~1$^{\prime }$ in
\cite{16}. The difference between the hypotheses occurs because of the right
continuity of the solution instead of left continuity.

\noindent More results on the oscillation of impulsive delay differential
equations can be found in the survey paper \cite{17}.

\noindent \textbf{Corollary 2. }Assume that $b(t)\equiv 0,~c(t)\neq 0$ and
that%
\begin{equation}
\lim \limits_{n\rightarrow \infty }\sup ~\int \limits_{n}^{n+1}\left( \prod
\limits_{n-1<t_{j}\leq s}(1-b_{j})\right) c(s)\exp \left( ~\int
\limits_{n-1}^{s}a(u)du\right) ds>1.  \label{12y}
\end{equation}%
Then every solution of Eq. (\ref{1})-(\ref{2}) is oscillatory.

\noindent \textbf{Remark 4.} If $b(t)\equiv 0,~c(t)\neq 0,$ then Eq. (\ref{1}%
)-(\ref{2}) reduces to an impulsive delay differential equation with
piecewise constant argument. Eq. (\ref{1})-(\ref{2}) with $b(t)\equiv 0,~$%
and $t_{i}=i,~i=1,2,...$ has been investigated in \cite{9}. So, Corollary 2
is a generalization of Theorem 4 in \cite{9}.

\noindent Moreover, in \cite{9}, a difference equation is a main tool for
the proofs. Similarly, in the other works such as \cite{4,5,6,7,10,11,18,19}
the relation between difference equations and differential equations with
piecewise constant arguments are underlined. Here, because of the existence
of continuous argument, we have diffuculty to obtain related difference
equation. So, we apply another technique which is worked for delay
differential equations.\medskip

\noindent \textbf{Theorem 3. }If one of the following conditions is
satisfied then every solution of equation (\ref{1})-(\ref{2}) is oscillatory:%
\begin{equation}
\lim \limits_{t\rightarrow \infty }\inf ~\int \limits_{t-l}^{t}\left( \prod
\limits_{s-\tau <t_{j}\leq s}(1-b_{j})\right) b(s)\exp \left( ~\int
\limits_{s-\tau }^{s}a(u)du\right) ds>\dfrac{1}{e},  \label{11}
\end{equation}%
\begin{equation}
\lim \limits_{n\rightarrow \infty }\inf ~\int \limits_{n+1-l}^{n+1}\left(
\prod \limits_{n-1<t_{j}\leq s}(1-b_{j})\right) c(s)\exp \left( ~\int
\limits_{n-1}^{s}a(u)du\right) ds>\dfrac{1}{e},  \label{11yy}
\end{equation}%
where $l=\min \{ \tau ,1\}.$\medskip

\noindent \textbf{Proof. }Let conditions (\ref{11}) or (\ref{11yy})\textbf{\
}is satisfied. We shall prove that the existence of eventually positive (or
negative) solutions leads to a contradiction. Let $x(t)$ be an eventually
positive solution of equation (\ref{1})-(\ref{2}). Then $y(t)=\prod%
\limits_{T<t_{j}\leq t}(1-b_{j})x(t)$ is an eventually positive solution of
equation (\ref{5}) such that $y(t)>0,~y(t-\tau )>0,~y([t-1])>0$ for $%
n+1>t\geq n>T.$ Using the same arguments in the proof of Theorem 2, we
obtained that $z(t)$ defined in (\ref{6y}) is nonincreasing for $t>T.$ We
consider two cases:

\textit{Case 1. }$\tau >1.$ Dividing inequality (\ref{8}) by $z(t),~$and
then integrating from $t-1$ to $t,$ it is obtained that%
\begin{equation}
\ln \frac{z(t-1)}{z(t)}\geq \int \limits_{t-1}^{t}P(s)\frac{z(s-1)}{z(s)}ds,
\label{12}
\end{equation}%
where $P(t)$ is defined in (\ref{6}). Since $e^{x}\geq ex$ for $x\in
\mathbb{R}
,$ we obtain that%
\begin{eqnarray}
\frac{z(t-1)}{z(t)} &\geq &\exp \left( ~\int \limits_{t-1}^{t}P(s)\frac{%
z(s-1)}{z(s)}ds\right)  \notag \\
&\geq &e\left( ~\int \limits_{t-1}^{t}P(s)\frac{z(s-1)}{z(s)}ds\right) .
\label{13}
\end{eqnarray}%
Let $u(t)=\dfrac{z(t-1)}{z(t)}.$ Since $z(t)$ is nonincreasing for $t>T$, $%
\underset{t\rightarrow \infty }{\lim \inf }~u(t)\geq 1.$

\noindent Assume that $\underset{t\rightarrow \infty }{\lim \inf }%
~u(t)=+\infty .$ Then integrating inequality (\ref{8}) from $t-\dfrac{1}{2}$
to $t,$ we have%
\begin{equation*}
z(t)-z(t-\frac{1}{2})+\int \limits_{t-\frac{1}{2}}^{t}P(s)z(s-1)ds\leq 0.
\end{equation*}%
Since $z(t)$ is nonincreasing, from the above inequality, we obtain that%
\begin{equation}
z(t)-z(t-\frac{1}{2})+z(t-1)\int \limits_{t-\frac{1}{2}}^{t}P(s)ds\leq 0.
\label{14}
\end{equation}%
Dividing inequality (\ref{14}) by $z(t)$ and $z(t-\frac{1}{2}),$ we get%
\begin{equation}
1-\frac{z(t-\frac{1}{2})}{z(t)}+\frac{z(t-1)}{z(t)}\int \limits_{t-\frac{1}{2%
}}^{t}P(s)ds\leq 0,  \label{15}
\end{equation}%
and%
\begin{equation}
\frac{z(t)}{z(t-\frac{1}{2})}-1+\frac{z(t-1)}{z(t-\frac{1}{2})}\int
\limits_{t-\frac{1}{2}}^{t}P(s)ds\leq 0,  \label{16}
\end{equation}%
respectively. Now, from (\ref{15}) we obtain%
\begin{equation*}
\underset{t\rightarrow \infty }{\lim \inf }\frac{z(t-\frac{1}{2})}{z(t)}%
=+\infty
\end{equation*}%
which contradicts with (\ref{16}). So,$\underset{t\rightarrow \infty }{\text{
}\lim \inf }~u(t)$ is finite.

\noindent If $\underset{t\rightarrow \infty }{\lim \inf ~}u(t)=w,~w\geq 1$
is finite, then inequality (\ref{13}) implies that%
\begin{equation*}
\lim \limits_{t\rightarrow \infty }\inf ~\int \limits_{t-1}^{t}P(s)ds\leq
\frac{1}{e}.
\end{equation*}%
In view of (\ref{6}), (\ref{8y}), and (\ref{9y}), we obtain from the above
inequality that%
\begin{equation*}
\lim \limits_{t\rightarrow \infty }\inf ~\int \limits_{t-1}^{t}\left( \prod
\limits_{s-\tau <t_{j}\leq s}(1-b_{j})\right) b(s)\exp \left( ~\int
\limits_{s-\tau }^{s}a(u)du\right) ds\leq \dfrac{1}{e},
\end{equation*}%
which contradicts the hypothesis (\ref{11}).

\noindent Now, dividing inequality (\ref{8}) by $z(t),~$and then integrating
from $n$ to $n+1,$ it is obtained that%
\begin{equation*}
\ln \frac{z(n)}{z(n+1)}\geq \int \limits_{n}^{n+1}P(s)\frac{z(s-1)}{z(s)}ds,
\end{equation*}%
where $P(t)$ is defined in (\ref{6}). Since $e^{x}\geq ex$ for $x\in
\mathbb{R}
,$ we obtain that%
\begin{equation}
\frac{z(n)}{z(n+1)}\geq e\left( ~\int \limits_{n}^{n+1}P(s)\frac{z(s-1)}{z(s)%
}ds\right) .  \label{18y}
\end{equation}%
Define $v(n)=\dfrac{z(n)}{z(n+1)}.$ Since $z(t)$ is nonincreasing for $t>T$,
$\underset{n\rightarrow \infty }{\lim \inf }~v(n)\geq 1.$ By doing the same
calculations with first part of the proof, we get that $\underset{%
n\rightarrow \infty }{\lim \inf }~v(n)$ is finite. Therefore, from the
inequality (\ref{18y}), we have
\begin{equation*}
\lim \limits_{n\rightarrow \infty }\inf ~\int \limits_{n}^{n+1}\left( \prod
\limits_{n-1<t_{j}\leq s}(1-b_{j})\right) c(s)\exp \left( ~\int
\limits_{n-1}^{s}a(u)du\right) ds\leq \dfrac{1}{e},
\end{equation*}%
which contradicts (\ref{11yy}).

\textit{Case 2.} $\tau \leq 1.$ Since the proof is similar to proof of
\textit{Case 1}, we shall give the sketch of the proof. Dividing inequality (%
\ref{9}) by $z(t),~$and then integrating from $t-\tau $ to $t,$ it is
obtained that%
\begin{equation}
\frac{z(t-\tau )}{z(t)}\geq e\left( ~\int \limits_{t-\tau }^{t}P(s)\frac{%
z(s-\tau )}{z(s)}ds\right) .  \label{19y}
\end{equation}%
Using the similar arguments in \textit{Case 1}, we get that $\underset{%
t\rightarrow \infty }{\lim \inf }~\dfrac{z(t-\tau )}{z(t)}$ is finite. So,
from inequality (\ref{19y}), we have%
\begin{equation*}
\lim \limits_{t\rightarrow \infty }\inf ~\int \limits_{t-\tau }^{t}\left(
\prod \limits_{s-\tau <t_{j}\leq s}(1-b_{j})\right) b(s)\exp \left( ~\int
\limits_{s-\tau }^{s}a(u)du\right) ds\leq \dfrac{1}{e},
\end{equation*}%
which contradicts (\ref{11}).

\noindent Moreover, dividing inequality (\ref{9}) by $z(t),~$and then
integrating from $n+1-\tau $ to $n+1,$ it is obtained that%
\begin{equation}
\frac{z(n+1-\tau )}{z(n+1)}\geq e\left( ~\int \limits_{n+1-\tau }^{n+1}P(s)%
\frac{z(s-\tau )}{z(s)}ds\right) .  \label{20y}
\end{equation}%
By using the similar arguments in \textit{Case 1}, we get that $\underset{%
n\rightarrow \infty }{\lim \inf }~\dfrac{z(n+1-\tau )}{z(t)}$ is finite. So,
from inequality (\ref{20y}), we have%
\begin{equation*}
\lim \limits_{n\rightarrow \infty }\inf ~\int \limits_{n+1-\tau
}^{n+1}\left( \prod \limits_{n-1<t_{j}\leq s}(1-b_{j})\right) c(s)\exp
\left( ~\int \limits_{n-1}^{s}a(u)du\right) ds\leq \dfrac{1}{e},
\end{equation*}%
which contradicts (\ref{11yy}). So, the proof is complete.

\noindent \textbf{Corollary 3.} Assume that $b(t)\neq 0$, $c(t)\equiv 0$ and
that%
\begin{equation*}
\lim \limits_{t\rightarrow \infty }\inf ~\int \limits_{t-\tau }^{t}\left(
\prod \limits_{s-\tau <t_{j}\leq s}(1-b_{j})\right) b(s)\exp \left( ~\int
\limits_{s-\tau }^{s}a(u)du\right) ds>\frac{1}{e}.
\end{equation*}%
Then every solution of Eq. (\ref{1})-(\ref{2}) is oscillatory.

\noindent \textbf{Corollary 4. }Assume that $b(t)\equiv 0,~c(t)\neq 0$ and
that%
\begin{equation*}
\lim \limits_{n\rightarrow \infty }\inf ~\int \limits_{n}^{n+1}\left( \prod
\limits_{n-1<t_{j}\leq s}(1-b_{j})\right) c(s)\exp \left( ~\int
\limits_{n-1}^{s}a(u)du\right) ds>\frac{1}{e}.
\end{equation*}%
Then every solution of Eq. (\ref{1})-(\ref{2}) is oscillatory.\medskip

\noindent Now, we give some examples to illustrate our results. Note that
previous results in the litarature can not be applied following differential
equations to obtain existence of oscillatory solutions.\smallskip

\noindent \textbf{Example 1. }Let us consider the following differential
equation%
\begin{equation}
\left \{
\begin{array}{c}
x^{\prime }(t)+\pi x(t-\dfrac{1}{2})+c(t)x([t-1])=0,~t\neq n,~n=1,2,...,~t>0,
\\
x(n^{+})-x(n^{-})=-x(n^{+}),~n=1,2,...,%
\end{array}%
\right.  \label{17}
\end{equation}%
where $c(t)\geq 0$ is any continuous function. It can be shown that the
hypotheses of Theorem 2 as well as Theorem 3 are satisfied. So, all
solutions of Eq. (\ref{17}) are oscillatory.

\noindent \textbf{Example 2. }Consider the following differential equation%
\begin{equation}
\left \{
\begin{array}{c}
x^{\prime }(t)+x(t)+\pi x(t-\dfrac{5}{2})+e^{t}x([t-1])=0,~t\neq
t_{n},~n=1,2,...,~t>0, \\
x(t_{n}^{+})-x(t_{n}^{-})=-2^{n}x(t_{n}^{+}),~n=1,2,...,%
\end{array}%
\right.  \label{18}
\end{equation}%
where $\{t_{n}\}_{n=1}^{\infty }$ is an increasing sequence such that $%
\underset{n\rightarrow \infty }{\lim }t_{n}=\infty .$

It is clear that $a(t)=1,~b(t)=\pi ,~c(t)=e^{t},~\tau =\dfrac{5}{2}$ and $%
b_{n}=-2^{n}.$ It can be shown that the hypotheses of Theorem 2 as well as
Theorem 3 are satisfied. So, all solutions of Eq. (\ref{18}) are oscillatory.

\end{document}